\documentclass[10pt,twoside]{siamart1116}

\usepackage[english]{babel}
\usepackage{graphicx,epstopdf,epsfig}
\usepackage{amsfonts,epsfig,fancyhdr,graphics, hyperref,amsmath,amssymb}

\usepackage{algorithm}%
\usepackage{algpseudocode}%
\usepackage{mathdots}

\newcommand{\HE}{Name of Handling Editor}
\newcommand{\DoS}{Month/Day/Year}
\newcommand{\DoA}{Month/Day/Year}
\newcommand{\CA}{Heike Fa\ss bender}

\newcommand{\Names}{Peter Benner, Heike Fa\ss bender, Michel-Niklas Senn}
\newcommand{\Title}{Structure-preserving Krylov Subspace Approximations for the Matrix Exponential of Hamiltonian Matrices: A Comparative Study}

\renewcommand{\theequation}{\thesection.\arabic{equation}}
\renewtheorem{theorem}{Theorem}[section]

\numberwithin{equation}{section} 

\begin{document}

\bibliographystyle{plain}

\setcounter{page}{1}

\thispagestyle{empty}

 \title{\Title\thanks{Received
 by the editors on \DoS.
 Accepted for publication on \DoA. 
 Handling Editor: \HE. Corresponding Author: \CA}}

\author{
Peter Benner\thanks{
Max Planck Institute for Dynamics of Complex Technical Systems,
39106 Magdeburg, Germany (benner@mpi-magdeburg.mpg.de)}
\and 
Heike Fa{\ss}bender\thanks{
Institute for Numerical Analysis, TU Braunschweig,  38106 Braunschweig, Germany (h.fassbender@tu-braunschweig.de)}
\and
Michel-Niklas Senn\thanks{
Institute for Numerical Analysis, TU Braunschweig, 38106 Braunschweig, Germany (michel.senn@outlook.de)}
}

\markboth{\Names}{\Title}

\maketitle

\begin{abstract}
We study structure-preserving Krylov subspace methods for approximating the matrix–vector products $f(H)b,$ where $H$ is a large Hamiltonian matrix and $f$ denotes either the matrix exponential or the related $\varphi$-function. Such computations are central to exponential integrators for Hamiltonian systems. Standard Krylov methods generally destroy the Hamiltonian structure under projection, motivating the use of Krylov bases with 
$J$-orthogonal columns that yield Hamiltonian projected matrices and symplectic reduced exponentials. We compare several such structure-preserving Krylov methods on representative Hamiltonian test problems, focusing on accuracy, efficiency, and structure preservation, and briefly discuss adaptive strategies for selecting the Krylov subspace dimension.
\end{abstract}


\section{Introduction}\label{sec1}
We present a comparative study of structure-preserving Krylov subspace methods for approximating the matrix–vector products
\begin{equation}\label{eq_problem}
  f(H)b = e^{H} b
  \quad \text{and} \quad
  f(H)b = \varphi(H) b,
\end{equation}
where
\[
  \varphi(z) = \frac{e^{z} - 1}{z}.
\]
The matrix $H\in \mathbb{R}^{2n \times 2n}$ is assumed to be \emph{Hamiltonian}, that is,
\[
  H J = (H J)^{T},
\]
for
\begin{equation}\label{eq_J}
J = J_n = \begin{bmatrix}0 & I_n\\-I_n&0\end{bmatrix}\in \mathbb{R}^{2n \times 2n},
\end{equation}
where $I_n$ is the $n \times n$ identity matrix. 
A key property in this setting is that the matrix exponential $e^{H}$
is \emph{symplectic} whenever $H$ is Hamiltonian, that is, $S^TJS=J$ holds for $S = e^H$.

Such problems arise, for example, in the numerical integration of Hamiltonian systems when exponential integrators are employed, where it is necessary to compute $f(H)b$ at each time step. In this context, $H$ is typically a large, sparse matrix obtained from the spatial discretization of an underlying partial differential equation describing Hamiltonian dynamics, while $f(H)$ is generally dense. Consequently, explicitly forming $f(H)$ is prohibitively expensive for high-dimensional problems. Therefore, we focus on approximating $f(H)b$ directly, without explicitly constructing the full matrix $f(H)$.

For general unstructured matrices, the problem of approximating the action of a matrix function has been investigated extensively in the literature, both as a stand-alone problem and in the context of exponential integrators; see, e.g., \cite{AH11,BerS24,GRT18,GoeG14,HLS98,HL97} and the references therein. Typically, the desired approximation is accomplished by using the well-known fact \cite[Theorem~1.13]{Hig08} that for $A \in \mathbb{R}^{\ell \times \ell}$ and any nonsingular matrix $X \in \mathbb{R}^{\ell \times \ell}$, $$f(X^{-1} A X) = X^{-1}f(A) X$$ holds,
provided that the function $f$ is defined on the spectrum of $A$. Assume that a suitable projection $\Pi = V W^{T} \in \mathbb{R}^{\ell \times \ell}$ is given, where $V, W \in \mathbb{R}^{\ell \times m}$, $m \leq \ell$, have full column rank and $W^{T} V = I_{m}$.
Then $f(A)b$ for $b \in \mathbb{R}^{\ell}$ can be approximated by
\begin{equation}\label{eq27a}
    f(A)b \approx V f(\widetilde{A}) W^{T} b
\end{equation}
with $\widetilde{A} = W^{T} A V \in \mathbb{R}^{m \times m}.$
If $V$ is chosen such that $V e_{1} = b / \| b \|_{2}$ (equivalently, $W^{T} b = \| b \|_{2} e_{1}$), this simplifies to
\begin{equation}\label{eq27}
    f(A)b \approx \| b \|_{2} V f(\widetilde{A}) e_{1}.
\end{equation}
When $m\ll \ell$, the projected matrix $\widetilde{A}$ is of much smaller dimension than the original matrix $A$; consequently, evaluating $f(\widetilde{A})b$ is far less computationally expensive than evaluating $f(A)b$. Thus, the function $f$ is evaluated on the small matrix argument $\widetilde{A} \in \mathbb{R}^{m \times m}$. The result is then mapped back to the original space $\mathbb{R}^{\ell}$. This approximation becomes exact when $\ell = m$. Typically, $V=W$ and the columns of $V$ form an orthonormal basis of  the standard \emph{Krylov subspace} $\mathcal{K}_m(A,b)=\text{span}\{b,Ab,\ldots,A^{m-1}b\}$ of order $m$. The matrix $V$ is usually computed via the standard Arnoldi algorithm, see, e.g., \cite{GolVL13,W07,W02} and the references therein. 

If $A = H \in \mathbb{R}^{2n \times 2n}$ is a Hamiltonian matrix, then the projected matrix $\widetilde H = V^THV \in \mathbb{R}^{m \times m}$ is in general not Hamiltonian. While an even value of $m$ is a necessary condition, it is not sufficient, $\widetilde H$ is typically non-Hamiltonian even then. In the context of symplectic exponential integrators (see, e.g., \cite{HLW06,MeiW17,KuoLS18,LR04,SC94,WuW19,EK19,LiC19}), it is often desirable to ensure that the projected matrix $\widetilde H = V^THW $ preserves the Hamiltonian structure such that $e^{\widetilde{H}}$ is symplectic. This can be achieved by choosing $m=2k$ even and $V$ as a matrix with $J$-orthogonal columns ($V^TJ_nV=J_k$) and $W^T= J_k^TV^TJ_n$ as its left inverse. 
 Hence, instead of the standard Arnoldi method, one should use a Krylov subspace method which generates a matrix $V_{2k}$ with $J$-orthogonal columns. 

Fortunately, there exist a number of generalizations of these standard Krylov subspace methods which aim at generating a  $J$-orthogonal basis. 
Building on the work of \cite{EK19, LiC19}, in which structure-preserving Krylov subspace methods for linear Hamiltonian systems were investigated, we consider the methods discussed there plus an additional one. However, our focus is solely on the approximation of $f(H)b$ as in \eqref{eq_problem} rather than their use in exponential integrators.

The remainder of this paper is structured as follows.
After providing some fundamental definitions and properties of Hamiltonian and symplectic matrices in Section \ref{sec2}, in Section \ref{sec3} we briefly introduce the Krylov subspace methods considered in this work and present their algorithmic realization explicitly for clarity and reproducibility.\footnote{In addition, our implementation is publicly available on Zenodo, see \cite{zenodo}.} Section \ref{sec:numerical_experiments_expint} describes six representative examples of Hamiltonian systems frequently encountered in the literature. For each example, the corresponding Hamiltonian matrix is given explicitly, so that the computation of the two functions in \eqref{eq_problem}, which arise in exponential integrators for such systems, can be reproduced. In Section \ref{sec5}, all Krylov subspace methods presented in Section \ref{sec3} are applied to the examples introduced in Section \ref{sec:numerical_experiments_expint} and compared with each other. Furthermore, we briefly demonstrate how the dimension of the Krylov subspace can be determined adaptively for the two most promising algorithms. We conclude with a summary and an outlook on potential further work in Section~\ref{sec6}.

\section{Fundamentals on Hamiltonian and symplectic matrices}\label{sec2}
It is easily checked that $J$ from \eqref{eq_J}  is orthogonal and skew-symmetric,
$J^T=J^{-1} =-J$. The matrix $J$ induces a skew-symmetric bilinear form $\langle\cdot,\cdot\rangle_J$ on $\mathbb{R}^{2n}$ defined by $\langle x,y \rangle_J = y^TJx$ for $x,y\in \mathbb{R}^{2n}.$

An equivalent definition for a matrix $H \in \mathbb{R}^{2n \times 2n}$ to be Hamiltonian is that there exist matrices $E$, $B=B^T$, $C=C^T$ $\in\mathbb{R}^{n\times n}$ such that
    \begin{align*}
        H = \begin{bmatrix}
            E & B \\
            C & -E^T
        \end{bmatrix}.
    \end{align*}
It is straightforward to check that Hamiltonian matrices are skew-adjoint with respect to the bilinear form $\langle\cdot,\cdot\rangle_J$, while symplectic matrices are orthogonal with respect to $\langle\cdot,\cdot\rangle_J$ and are therefore also called \emph{J-orthogonal}. The $2n \times 2n$ symplectic matrices form a Lie group, the $2n \times 2n$ Hamiltonian matrices the associated Lie algebra. Numerous properties of the sets of Hamiltonian and symplectic matrices (and their interplay) have been studied in the literature, see, e.g., \cite{MMT03} and the references therein. Here we list only those properties necessary for the following discussion.

\begin{lemma}\label{lem1}
\begin{enumerate}
    \item Let $H$ be a nonsingular Hamiltonian matrix. Then $H^{-1}$ is Hamiltonian as well.
    \item Let $S$ be a symplectic matrix. Then $S^{-1} = -JS^TJ$ is symplectic as well.
    \item Let $H$ be a Hamiltonian matrix and $S$ be a  symplectic matrix. Then $S^{-1}HS$ is a Hamiltonian matrix.
    \item 
    An orthogonal matrix $S\in\mathbb{R}^{2n\times2n}$ is symplectic if and only if
    \begin{align}\label{eq:Unitary_Symplectic_Matrices}
        S = \begin{bmatrix}
            V_1 & V_2\\
            -V_2 & V_1
        \end{bmatrix},
    \end{align}
    where $V_1,V_2\in\mathbb{R}^{n\times n}$ with $V_1^TV_1+V_2^TV_2 =I_n$ and $V_1^TV_2 = V_2^TV_1$, see, e.g, \cite{PV81}.
    \end{enumerate}
\end{lemma}
In the context of Krylov subspace methods, we will usually do not generate a full (orthogonal and/or symplectic) basis for $\mathbb{R}^{2n}$. Thus, we will consider matrices $S \in \mathbb{R}^{2n \times 2k}$, $k \leq n$, with $J$-orthogonal columns. Any matrix $S \in \mathbb{R}^{2n \times 2k}$, $k \leq n$ with $S^TJ_nS=J_k$ is said to  have \emph{J-orthogonal columns}. The following statements can be verified in a straightforward way.
\begin{lemma}\label{lem3}
Let $S \in \mathbb{R}^{2n \times 2k}$, $k \leq n$, be a matrix with J-orthogonal columns, $S^TJ_nS=J_k.$ 
Let $H \in \mathbb{R}^{2n \times 2n}$ be Hamiltonian.
\begin{enumerate}
    \item The matrix $J_k^TS^TJ_n\in \mathbb{R}^{2k \times 2n}$ is the left inverse of $S$, $J_k^TS^TJ_nS = I_{2k}.$
    \item The matrix $(J_k^TS^TJ_n)HS \in \mathbb{R}^{2k \times 2k}$ is Hamiltonian.
\end{enumerate}
\end{lemma}

\section{Krylov subspace projection methods}\label{sec3}
A standard approach to approximate $f(H)b$ employs a matrix $U_{2k} \in \mathbb{R}^{2n \times 2k}$ ($k \le n$ chosen appropriately)
\begin{align}\label{eq:FuncApprox}
    f(H) b
    \;\approx\;
     U_{2k} f(H_{2k}) W_{2k}^{T} b,
\end{align}
see \cite{GS92},
where $H_{2k} = W_{2k}^TAU_{2k}$ and $W_{2k}^TU_{2k}=I_{2k},$ that is, e.g., 
\[
   W_{2k}^{T} =
   \begin{cases}
      U_{2k}^{T}, & \text{if } U_{2k} \text{ is orthogonal}, \\
      J_k^{T} U_{2k}^{T} J_n, & \text{if } U_{2k} \text{is } J\text{-ortogonal}.
   \end{cases}
\]

 The standard Krylov subspace $\mathcal{K}_m(H,b)=\text{span}\{b,Hb,\ldots,H^{m-1}b\}$ of order $m$ comprises all polynomials of $H$ up to degree $m-1$, multiplied by $b,$ that is, for all polynomials $p_{m-1}$ of degree $\leq m-1$ we have
 \[
 p_{m-1}(H)b = V_mp_{m-1}(W_m^THV_m)e_1
 \]
 in case the columns of $V_m$ span $\mathcal{K}_m(H,b),$ $V_me_1 = b/\|b\|_2,$ and $W_m^TV_m=I_m.$ Typically, an orthonormal basis of  $\mathcal{K}_m(H,b)$  is computed via the standard Arnoldi algorithm, see, e.g., \cite{GolVL13,W07,W02} and the references therein. As such a  Krylov subspace approximation to  $f(H)b$ usually converges fast, see, e.g., \cite{Hig08,HL97} and the references therein,  a quintessential choice of $V$ and $W$ in \eqref{eq27} is the orthonormal matrix $V_m=V=W$ generated by the Arnoldi method for the Krylov subspace $\mathcal{K}_m(H,b)$. 
But in this case the projected matrix $\widetilde H = V_m^THV_m \in \mathbb{R}^{m \times m}$ is generally not Hamiltonian. Structure preservation is enforced by requiring that $\widetilde H = W^THV \in \mathbb{R}^{2k \times 2k}$ is a Hamiltonian matrix. This can be achieved by choosing $V=S\in \mathbb{R}^{2n \times 2k}$ as a matrix with J-orthogonal columns and $W^T= J_k^TS^TJ_n$ as its left inverse due to Lemma \ref{lem3}. Hence, instead of the standard Arnoldi method one should use a Krylov subspace method which generates a matrix $V$ with $J$-orthogonal columns. 
Fortunately, there exist a number of generalizations of these standard Krylov subspace methods which aim at generating a  $J$-orthogonal basis.
In the following, we list a number of such methods and briefly state some of their features. We refer to the existing literature for all further details.  For clarity and reproducibility, we provide listings of all corresponding algorithms in the following subsections.

Block Krylov subspace methods for approximating $e^{H}B$, where $H \in \mathbb{R}^{2n \times 2n}$ is a Hamiltonian and skew-symmetric matrix and $B \in \mathbb{R}^{2n \times \ell}$ has only a few columns, have been studied in \cite{ArcBA20,ArcB22}. In \cite{ArcBA20}, a structure-preserving block Lanczos process is employed to reduce the Hamiltonian matrix to the form $\left[
\begin{smallmatrix} 
A & N \\ 
M & -A^T 
\end{smallmatrix} \right],$
where $A, N, M$ are block tridiagonal matrices of size $\ell k \times \ell k$ with $\ell \times \ell$ blocks, and $N$ and $M$ 
are additionally symmetric. In contrast, \cite{ArcB22} applies a global Hamiltonian Lanczos process, leading to a reduced Hamiltonian matrix of the same block form,  but with $A \in \mathbb{R}^{k \times k}$ tridiagonal and $N, M \in \mathbb{R}^{k \times k}$ diagonal. Another idea has been proposed in \cite{LopS06}. In the present study, we consider only the computation of  $f(H)b$ for a single vector $b$, and not 
$f(H)B$ for multiple vectors.

\subsection{Arnoldi Method (A)}
Before discussing the structure-preserving methods, we briefly recall the standard approach, namely the Arnoldi method, since we will subsequently compare the structure-preserving schemes with this classical technique.
The method can be applied to any square matrix; therefore, in this section, we consider $H \in \mathbb{R}^{n \times n}$. The Arnoldi method generates an orthonormal basis of the standard \emph{Krylov subspace} of order $k \in \mathbb{N}$ associated with $H \in\mathbb{R}^{n \times n}$ and a given vector $b\in\mathbb{R}^{n}\setminus\{0\}$, given by
\[
\mathcal{K}_{k}(H,b):=\text{span}\left\{b,Hb,H^2b,\dots,H^{k-1}b\right\},
\]
see, e.g., \cite{GolVL13,W02,W07} and the references therein. In case the method does not stop prematurely, after $k$ steps it has computed an $n \times k$ matrix $U_{k}$ with orthonormal columns such that
\begin{equation}
HU_{k} = U_{k}\widehat{H}_{k} + u_{k+1}\widehat{h}_{k+1,k}e_{k}^T = U_{k+1} \widetilde{H}_{k+1} \label{eq:arnoldi_recursion}
\end{equation}
with a $k \times k$ upper Hessenberg matrix $\widehat{H}_{k} = U_{k}^THU_{k}$, $U_{k}e_1 = b/\|b\|_2$ and $\widetilde{H}_{k+1} = \left[ \begin{smallmatrix} \widehat{H}_{k} \\ \widehat{h}_{k+1,k}e_{k}^T\end{smallmatrix} \right] \in \mathbb{R}^{k+1 \times k}$. The column $u_j$ of $U_{k}$ can be computed with the help of all preceding columns $u_1, \ldots, u_{j-1}.$ Algorithm \ref{alg:A} summarizes the Arnoldi process.
 As is well known, the algorithm (as all algorithms considered here) may fail to compute a basis of the desired length due to the condition on Line 9, which checks whether $\vert \widehat{h}_{j+1,j}\vert$ is approximately zero using the user-defined tolerance \texttt{tol}. Note that if the Arnoldi method stops in Line~10, then an (approximate) $H$-invariant subspace is found. In this case, $b$ is contained in an $H$-invariant subspace, so that the approximation from \eqref{eq27} yields the exact (up to numerical rounding errors) result in \eqref{eq_problem}. 
 For the examples considered here, and for the relatively small subspace dimensions chosen, this issue never arose. It is also well known that the computed vectors gradually lose their (theoretically exact) orthogonality as the iteration proceeds, so that reorthogonalization may become necessary.

\begin{algorithm}[ht!]
    \caption{Arnoldi Process}
    \label{alg:A}
    \begin{algorithmic}[1]
        \Require $H\in\mathbb{R}^{n\times n}$,  $v\in\mathbb{R}^{n}$, $k\in\mathbb{N}$, $\textrm{tol}\in\mathbb{R}_{>0}$
        \Ensure  $\widetilde{H}_{k+1}\in\mathbb{R}^{(k+1)\times k}$, $U_{k+1}\in\mathbb{R}^{n\times (k+1)}$ with $U_{k+1}^TU_{k+1} =I_{k+1}$
        
        \State $u_{1} = v/\Vert v\Vert_2$
        \For {$j=1,2,\dots,k$}
        \State $z=Au_{j}$
        \For{$i=1,2,\dots,j$}
        \State $\widehat{h}_{i,j}=u_{i}^Tz$
        \State $z = z - u_{i}\widehat{h}_{i,j}$
        \EndFor
        \State $\widehat{h}_{j+1,j} = \Vert z\Vert_2$
        \If{$\widehat{h}_{j+1,j}<\textrm{tol}$}
        \State stop
        \EndIf
        \State $u_{j+1} = z/\widehat{h}_{j+1,j}$
        \EndFor
    \end{algorithmic}
\end{algorithm}

In case $H$ is symmetric, the matrix $\widehat{H}_k$ is symmetric and thus tridiagonal. The computation of $u_j$ simplifies, only the columns $u_{j-1}$ and $u_{j-2}$ are required ('three-term-recurrence'). The resulting method is called \emph{Lanczos method}. This method can be extended for unsymmetric matrices to the \emph{unsymmetric Lanczos method} which computes matrices $Q$ and $P$ via a short recurrence so that $P^THQ=T$ is tridiagonal and $P^TQ=I_n$, see, e.g. \cite{GolVL13,W07}. 

Since we will be applying the algorithm to Hamiltonian matrices $H \in \mathbb{R}^{2n \times 2n}$, we always let the algorithm run for an even number of steps, i.e., $k = 2r$. In this case, the Arnoldi method requires $2r$ matrix-vector-multiplications and $r^2$ inner products to generate the orthonormal basis.
Recall that even if $H$ is Hamiltonian, $\widehat H_{2r}$ and $\widehat{T}_{2r}$ will in general not be Hamiltonian.

\subsection{Hamiltonian Lanczos Method (HL)}
Applied to a  Hamiltonian matrix $H \in \mathbb{R}^{2n \times 2n}$, this method (see Algorithm \ref{alg:HL} for a summary) generates a matrix $S_{2k} = \left[U_k ~~V_k\right],$ $U_k, V_k \in \mathbb{R}^{2n \times k}$ with $J$-orthogonal columns such that
\begin{equation} \label{eq:hamlan_recursion}
H[U_k~V_k] = [U_k~V_k] \begin{bmatrix} G_k & T_k\\ D_k & -G_k\end{bmatrix} + u_{k+1}\beta_{k}e_{2k}^T,
\end{equation}
where $ G_k = \operatorname{diag}(\gamma_1, \ldots, \gamma_k), D_k = \operatorname{diag}(\delta_1, \ldots, \delta_k) \in \mathbb{R}^{k \times k}$ are diagonal matrices and 
\[
T_k = \left[\begin{smallmatrix}\alpha_1 & \beta_1 &\\ \beta_1 & \ddots & \ddots\\ &\ddots & \ddots & \beta_{k-1}\\ &&\beta_{k-1} & \alpha_k \end{smallmatrix}\right] \in \mathbb{R}^{k \times k}
\]
is a symmetric tridiagonal matrix.
 Like the (un)symmetric Lanczos method, only a short recurrence is needed to compute the next vectors $u_{k+1}$ and $v_{k+1}$ of the basis involving only the three preceding vectors $v_k, u_k, u_{k-1}.$ The columns of $S_{2k}$ form a symplectic basis of the standard Krylov subspace $\mathcal{K}_{2k}(H,u_1).$ 
There is some degree of freedom in the algorithm that we exploit here to minimize the condition number of $S_{2k}$, as suggested in \cite{BFS11}. To enable a fair comparison with the other methods, we choose $k=r$, such that a Krylov subspace of dimension  $2r$ is generated (provided that no breakdown occurs). The method then requires $2r$ matrix-vector products and $3r$ inner products. See \cite{BenF97,BFS11,W04} for more details and different variants of the method.

\begin{algorithm}[ht!]
    \caption{Hamiltonian Lanczos Process}
    \label{alg:HL}
    \begin{algorithmic}[1]
        \Require Hamiltonian $H{\in}\mathbb{R}^{2n\times 2n}$, $\widehat{u}_1{\in}\mathbb{R}^{2n}$, $k{\in}\mathbb{N}$, $\textrm{tol}\in\mathbb{R}_{>0}$
        \Ensure  Hamiltonian $H_{2k} {\in}\mathbb{R}^{2k\times 2k}$, $S_{2k}{\in}\mathbb{R}^{2n\times 2k}$ with $S_{2k}^TJ_nS_{2k}=J_{k}$
        \State $u_0 = 0$
        \State $\beta_0 = \Vert \widehat{u}_1 \Vert_2$
        \State $u_1 = \widehat{u}_1/\beta_{0}$
        \For {$j=1,2,\dots,k$}
        \State $u = Hu_{j}$
        \State $\gamma_j = u_{j}^Tu$
        \State $v_{j} = u - \gamma_ju_{j}$
        \State $\delta_j = u^TJ_nu_{j}$
        \If{$|\delta_j|<$ tol}
        \State stop
        \EndIf
        \State $v_{j} = v_{j}/\delta_j$
        \State $v = Hv_j$
        \State $\alpha_j = -v^TJ_nv_{j}$
        \State $u_{j+1} = v - \beta_{j-1}u_{j-1} - \alpha_ju_{j} + \gamma_jv_{j}$
        \State $\beta_{j} = \Vert u_{j+1}\Vert_2$
        \If{$\beta_j<$ tol}
        \State stop
        \EndIf
        \State $u_{j+1} = u_{j+1}/\beta_{j}$
        \EndFor
        \State $S_{2k} = \begin{bmatrix}U_{k} & V_k\end{bmatrix}$
        \State $H_{2k} = \left[ \begin{smallmatrix} G_k & T_k\\ D_k & -G_k\end{smallmatrix}\right]$
    \end{algorithmic}
\end{algorithm}

\subsection{Symplectic Arnoldi Method (SA)}
In \cite{EK19}, a simple modification of the standard Arnoldi method for Hamiltonian matrices $H$ is suggested. Besides the orthogonal matrix $U_k$, a second matrix $V_k$ is generated by reorthogonalizing each column of $U_k$ against all vectors in $V_k$ with respect to $\langle \cdot, \cdot \rangle$ and $\langle \cdot, \cdot \rangle_J.$
The columns of the matrix $S_{2k} = [V_k ~-J_nV_k]$ form an orthogonal symplectic basis of $\mathcal{K}_{k}(H,v_1)+J_n\mathcal{K}_k(H,v_1),$ while the columns of $U_k$ are still an orthonormal basis of $\mathcal{K}_{k}(H,u_1).$ Here, $u_1 =U_k e_1 = v_1 =V_ke_1$ holds. Moreover, $S_{2k}^THS_{2k} = H_{2k}$ is a $2k \times 2k$ Hamiltonian matrix with no additional structure. 
 To enable comparison with the other algorithms, we formally regard $\mathcal{K}_{k}(H,v_1)+J_n\mathcal{K}_k(H,v_1)$ as a subspace having dimension $2k$. Therefore, we choose $k = r$.
For this choice, the method requires $r$ matrix-vector products and $r^2$ inner products to generate $V_r$. In addition, $S_{2r}^THS_{2r}$ needs to be computed. Our version of the algorithm is summarized in Algorithm \ref{alg:SA}.

\begin{algorithm}[ht!]
    \caption{Symplectic Arnoldi Process}
    \label{alg:SA}
    \begin{algorithmic}[1]
        \Require Hamiltonian $H{\in}\mathbb{R}^{2n\times2n}$, $v{\in}\mathbb{R}^{2n}$, $k{\in}\mathbb{N}$, $\textrm{tol}\in\mathbb{R}_{>0}$
        \Ensure  Hamiltonian $H_{2k}{\in}\mathbb{R}^{2k\times2k}$, $S_{2k}{\in}\mathbb{R}^{2n\times2k}$ with $S_{2k}^TJ_nS_{2k} =J_k$, $S_{2k}^TS_{2k}=I_{2k}$
        \State $u_1 = v/\Vert v\Vert_2$
        \State $v_1 = u_1$
        \For {$j=1,2,\dots,k$}
        \State $z = Hu_{j}$
        \For{$i = 1,2,\dots,j$}
            \State $z = z - (u_i^Tz) u_i$
        \EndFor
        \If{$\Vert z\Vert_2 <$ tol}
            \State stop
        \EndIf
        \State $u_{j+1} = z/\Vert z\Vert_2$
        \State $v = u_{j+1}$
        \For{$i = 1,2,\dots,j$}
            \State $v = v - (v_i^Tv) v_i - (v_i^TJ_n^Tv) J_nv_i$
            \State $a_{i,j}= v_i^THv_j$
            \State $a_{j,i}= v_j^THv_i$
            \State $m_{i,j}=m_{j,i}= -v_i^TJ_n^THv_j$
            \State $n_{i,j}=n_{j,i}= -v_i^THJ_nv_j$
        \EndFor
        \If{$\Vert v\Vert_2 <$ tol}
            \State stop
        \EndIf
        \State $v_{j+1} = v/\Vert v\Vert_2$
        \EndFor
        \State $S_{2r}=\begin{bmatrix}
            V_r & {-}J_nV_r
        \end{bmatrix}$
        \State $H_{2k}= \left[\begin{smallmatrix}A_k &N_k \\M_k & -A_k\end{smallmatrix}\right]$ 
    \end{algorithmic}
\end{algorithm}

\subsection{Isotropic Arnoldi Method (IA)}\label{subsec3}
 In \cite{MW00}, it is suggested to modify the standard Arnoldi algorithm by supplementing the orthogonalization of the vectors with an additional $J$-orthogonalization.
 This yields a matrix $S_{2k} = [U_k ~-J_nU_k]$, $U_k \in \mathbb{R}^{2n \times k}, U_k^TU_k = I_k$, with orthogonal and $J$-orthogonal columns such that
 \[
   J_k^TS_{2k}^TJ_n HS_{2k} = \begin{bmatrix}
        U_k & -J_nU_k
    \end{bmatrix}\begin{bmatrix}
        {T}_k & N_k\\
        -D_k & -T_k^T
    \end{bmatrix} = H_{2k},
\]
where $T_k \in\mathbb{R}^{k\times k}$ is upper Hessenberg, $D_k  \in\mathbb{R}^{k\times k}$ is diagonal and $N_k=N_k^T.$ As with the standard Arnoldi method, the calculation of the following vector, $u_{k+1}$, requires all previously computed vectors $u_1, \ldots, u_k$, and the vectors, $Ju_1, \ldots, Ju_k$. 

The columns of the orthonormal matrix $U_k$ span an isotropic subspace, hence the name isotropic Arnoldi process. The columns of $S_{2k}$ in general neither contain the Krylov subspace $\mathcal{K}_{k}(H,u_1)$ nor span another Krylov subspace. Nevertheless, we include this method in our study, since it is based on the Krylov subspace $\mathcal{K}_{k}(H,u_1)$. Although this classification is formally not correct, we treat it as a Krylov subspace method and, for comparison with the other approaches, choose $k=r$. This approach necessitates $r$ matrix-vector products and $r^2$ inner products.

\begin{algorithm}[ht!]
    \caption{Isotropic Arnoldi Process}
    \label{alg:IA}
    \begin{algorithmic}[1]
        \Require Hamiltonian $H{\in}\mathbb{R}^{2n\times2n}$, $v{\in}\mathbb{R}^{2n}$, $k{\in}\mathbb{N}$, $\textrm{tol}\in\mathbb{R}_{>0}$
        \Ensure  Hamiltonian $H_{2k}{\in}\mathbb{R}^{2k\times2k}$, $S_{2k}{\in}\mathbb{R}^{2n\times2k}$ with $S_{2k}^TJ_nS_{2k} =J_k$, $S_{2k}^TS_{2k}=I_{2k}$
        \State $u_1 = v/\Vert v\Vert_2$
        \For {$j=1,2,\dots,k$}
            \State $z = Hu_ {j}$
            \For{$i=1,2,\dots,j$}
                \State $t_{i,j} = u_i^Tz$
                \State $z = z - u_it_{i,j}$
            \EndFor
            \State $d_{j,j} = u_j^TJ_n^Tz$
            \State $z = z - J_nu_jd_{j,j}$
            \State $t_{j+1,j}= \Vert z\Vert_2$
            \If{$t_{j+1,1}<$ tol}
                \State stop
            \EndIf
            \State $u_{j+1} = z/t_{j+1,j}$
            \State $z = HJ_nu_{j}$
            \State $n_{j,j} = -u_j^Tz$
            \For{$i=1,2,\dots,j-1$}
                \State $n_{i,j} = n_{j,i} = -u_i^Tz$
            \EndFor
        \EndFor
        \State $S_{2k}=\begin{bmatrix}
            U_k & -J_nU_k
        \end{bmatrix}$
        \State $H_{2k}=\left[\begin{smallmatrix}
            T_k & N_k\\ -D_k& T_k
        \end{smallmatrix}\right]$ 
    \end{algorithmic}
\end{algorithm}

\subsection{Hamiltonian Extended Krylov Subspace Method (HEKS)}\label{subsec4}
In \cite{M11} it is suggested to consider the extended Krylov subspace
\[
    \mathcal{E}\mathcal{K}_{2k}(H,u) =\mathcal{K}_{2t}(H,u)+\mathcal{K}_{2s}(H^{{-}1},H^{{-}1}u),
\]
where $2k = 2(t+s)$ and either $t=s$ or $t=s+1$.
In case $H$ is a nonsingular Hamiltonian matrix, $H^{-1}$ is Hamiltonian as well (Lemma \ref{lem1}). Now, a basis  $\left[ U_t ~~ V_t\right]$ for $\mathcal{K}_{2t}(H,u)$ and another basis $\left[X_s ~~ Y_s \right]$ for $\mathcal{K}_{2s}(H^{-1},H^{-1}u)$ are generated such that each new vector is $J$-orthogonalized against all previously generated vectors. The columns of the matrix $S_{2k} = \left[ Y_s~~ U_t~~X_s~~V_t\right] \in \mathbb{R}^{2n \times 2k}$ are $J$-orthogonal  such that
\begin{equation}\label{eq:struct1}
 J_{k}^TS_{2k}^{T}J_nHS_{2k} = \begin{bmatrix}
0 & 0 & \Lambda_{s} & B_{st}\\
0 & 0 & B_{st}^T & T_{t}\\
\Delta_{s} & 0 & 0 & 0\\
0 & \Theta_{t} &0 & 0
\end{bmatrix}\in \mathbb{R}^{2k \times 2k}
\end{equation}
with diagonal matrices
\begin{align*}
    \Delta _s &= \operatorname{diag}(\delta_s, \ldots, \delta_1) \in \mathbb{R}^{s \times s}, \\
    \Theta_t &=\operatorname{diag}(\vartheta_1, \ldots, \vartheta_t) \in \mathbb{R}^{t \times t}, \\
    \Lambda_s &=\operatorname{diag}(\lambda_s, \ldots, \lambda_1) \in \mathbb{R}^{s \times s}, 
\end{align*}
a symmetric tridiagonal matrix
\begin{align*}
    T_t &= \left[\begin{smallmatrix}
        \alpha_1 &\beta_2&& \\
        \beta_2 & \ddots &\ddots& \\
        & \ddots &\ddots& \beta_t\\
        &  &\beta_t& \alpha_t
    \end{smallmatrix} \right]\in \mathbb{R}^{t \times t}, 
\end{align*}
and an anti-lower bidiagonal matrix $B_{st} \in \mathbb{R}^{s \times t}$ (depending on the choice of $t$)
\begin{align*}
B_{t-1,t} =\left[\begin{smallmatrix}
    &&&\gamma_{t-1} & \mu_{t}\\
    &  &\iddots &  \mu_{t-1}& \\
    & \iddots & \iddots &  & \\
    \gamma_{1} & \mu_{2} &  & &
\end{smallmatrix}\right]\in \mathbb{R}^{(t-1)\times t} 
\quad \textrm{or}\quad
B_{tt} =\left[\begin{smallmatrix}
    &&&\gamma_{t} \\
    &  &\iddots &  \mu_{t} \\
    & \iddots & \iddots &   \\
    \gamma_{1} & \mu_{2} &
\end{smallmatrix}\right]\in \mathbb{R}^{t\times t}, 
\end{align*}
see \cite{BFS22} for details. It is important to note that $S_{2k}e_1 \neq u~ (=S_{2k}e_{s+1}).$

For the choice $t=s=\ell$ (that is, $k = 2\ell$), the method requires $4\ell=2k$ matrix-vector products with $H$, $3\ell-1 = \frac{3}{2}k-1$ solves of linear systems with $H$ (efficiently implemented in the form $(JH)x = Jb$ making use of the symmetry of $JH$) and $14\ell=7k$ inner products. Evidently, $k=r=2\ell$ is chosen for the comparison of the algorithms, since the columns of  $S_{2r}$ span the Krylov subspace $ \mathcal{E}\mathcal{K}_{2r}(H,u)$.

{\small
    \begin{algorithm}[ht!]
        \caption{HEKS}
        \label{alg:heks}
        \begin{algorithmic}[1]
            \Require Hamiltonian matrix $H{\in}\mathbb{R}^{2n\times2n}$, $u_1{\in}\mathbb{R}^{2n}$ with $\|u_1\|_2=1$, $k{\in}\mathbb{N}$
            \Ensure for $t=s=\ell$:  $S_{2r} = S_{4\ell}= [y_\ell\;\dots\;y_1\;u_1\;\dots\;u_\ell\; x_\ell\;\dots\;x_1\;v_1\;\dots\;v_\ell]\in \mathbb{R}^{2n \times 4\ell}$ with $S_{4\ell}^TJ_nS_{4\ell}=J_{2\ell}$ and $H_{4\ell} =  J_{2\ell}S_{4\ell}^TJ_nHS_{4\ell}$ as in \eqref{eq:struct1}\newline
            (for $t=s-1=\ell$: $S_{2(\ell+1)}\in \mathbb{R}^{2n \times (4\ell+2)}$ the algorithm needs to be modified appropriately)
            
            \State $u_1 = u_1/\|u_1\|_2$  \Comment{Set up $S_1 = [u_1\mid v_1]$}
            \State $\vartheta_1 = u_1^TJ_nHu_1$
            \State $v_1 = Hu_1/\vartheta_1$
            \State $f_{11} = u_1^TJ_nH^{-1}u_1$\Comment{Set up $S_2 = [y_1~u_1\mid x_1~v_1]$}
            \State $w_x = H^{-1}u_1 -f_{11}v_1$
            \State $x_1 = w_x/\|w_x\|_2$
            \State $y_1 = H^{-1}x_1/\big(x_1^TJ_nH^{-1}x_1\big)$
            \State $\lambda_1 = -x_1^TJ_nHx_1$ and  $\delta_1 = y_1^TJ_nHy_1$
            \State $\alpha_1 = -v_1^TJ_nHv_1$ and  $\gamma_1 = -x_1^TJ_nHv_1$\Comment{Set up $S_3 = [y_1~u_1~u_2\mid x_1~v_1~v_2]$}
            \State $w_u = Hv_1-\gamma_1y_1-\alpha_1u_1$
            \State $u_2 = w_u/\|w_u\|_2$
            \State $\vartheta_2 = u_2^TJ_nHu_2$
            \State $v_2 = Hu_2/\vartheta_2$
            \State $e_{11} = y_1^TJ_nH^{-1}y_1$\Comment{Set up $S_4 = [y_2~y_1~u_1~u_2\mid x_2~x_1~v_1~v_2]$}
            \State $ g_{11} = y_1^TJ_nH^{-1}u_1,$ and $g_{12}=y_1^TJ_nH^{-1}u_2$
            \State $w_x = H^{-1}y_1 -e_{11}x_1 - g_{11}v_1-g_{12}v_2$
            \State $x_2 = w_x/\|w_x\|_2$
            \State $y_2 = H^{-1}x_2/\big((H^{-1}x_2)^TJ_nx_2\big)$
            \State $\lambda_2 = -x_2^TJ_nHx_2$ and $\delta_2 = y_2^TJ_nHy_2$
            \For{$j=3, 4, \ldots, \ell$}
            \State $\alpha_{j-1} = -v_{j-1}^TJ_nHv_{j-1}$ and  $\beta_{j-1} = -v_{j-1}^TJ_nHv_{j-2}$\Comment{Set up $S_{2j-1}$}
            \State $\gamma_{j-1}= -x_{j-1}^TJ_nHv_{j-1}$ and $ \mu_{j-1}=-x_{j-2}^TJ_nHv_{j-1}$
            \State $w_u = Hv_{j-1}-\gamma_{j-1} y_{j-1}-\mu_{j-1}y_{j-2}-\beta_{j-1} u_{j-2}-\alpha_{j-1} u_{j-1}$
            \State $u_j = w_u/ \|w_u\|_2$
            \State $\vartheta_j = u_j^TJ_nHu_j$
            \State $v_j = Hu_j/\vartheta_j$
            \State $g_{j-1,j-1}= y_{j-1}^TJ_nH^{-1}u_{j-1}$ and $g_{j-1,j}=y_{j-1}^TJ_nH^{-1}u_{j}$ \Comment{Set up $S_{2j}$}
            \State $e_{j-1,j-1}= y_{j-1}^TJ_nH^{-1}y_{j-1}$ and $e_{j-2,j-1}=y_{j-1}^TJ_nH^{-1}y_{j-2}$
            \State $w_x = H^{-1}y_{j-1} -e_{j-1,j-1}x_{j-1} -e_{j-2,j-1}x_{j-2}- g_{j-1,j-1}v_{j-1}-g_{j-1,j}v_j$
            \State $x_j = w_x/\|w_x\|_2$
            \State $y_j = H^{-1}x_j/\big((H^{-1}x_j)^TJ_nx_j\big)$
            \State $\lambda_j = -x_j^TJ_nHx_j$ and $\delta_j = y_j^TJ_nHy_j$
            \EndFor
            \State $\alpha_\ell = -v_\ell^TJ_nHv_\ell$ and $\beta_\ell = -v_\ell^TJ_nHv_{\ell-1}$
            \State $\gamma_\ell =-x_\ell^TJ_nHv_\ell$ and $ \mu_\ell=-x_{\ell-1}^TJ_nHv_\ell$
        \end{algorithmic}
    \end{algorithm}
}

\subsection{Block $J$-orthogonal Basis Method (BJ)}\label{subsec5}
In \cite{LiC19}, the standard Arnoldi method is used to generate an orthonormal matrix $U_k \in \mathbb{R}^{2n \times k}$ whose columns form a basis for $\mathcal{K}_{k}(H,v).$ Then $U_k$ is split into its first $n$ and last $n$ rows, $U_k = \left[\begin{smallmatrix} U_k^u\\ U_k^l\end{smallmatrix}\right],$ $U_k^u,$ $U_k^l \in \mathbb{R}^{n \times k}$ and rearranged as $V_{2k} = [U_k^u~~U_k^l] \in \mathbb{R}^{n \times 2k}.$ Next, the columns of $V_{2k}$ are orthogonalized to obtain a matrix $W_{2k}$ with orthonormal columns, $W_{2k}^TW_{2k}=I_{2k}$. For ease of notation, we assume that all columns of $[U_k^u~~U_k^l]$ are linearly independent, such that $W_{2k}$ has $2k$ columns. But in practice, $W_{2k}$ may have just $2p$ columns with $p \leq k$. Then
\[
S_{4k} = \begin{bmatrix} W_{2k} & 0\\ 0 & W_{2k}\end{bmatrix}\in \mathbb{R}^{2n \times 4k}
\]
has $J$-orthogonal columns, see Algorithm \ref{alg:BJ} for details. It is important to note that $S_{4r}e_1 \neq v.$

The Krylov subspace $\mathcal{K}_{k}(H,v)$ is contained in the range of $S_{4k}$.
To ensure comparability, and because all methods under consideration ideally construct a Krylov subspace of dimension $2r$, we set $k=2r$ in this case. Then, the Arnoldi method requires $2r$ matrix-vector-multiplications and $r^2$ inner products to generate the orthonormal basis $U_r.$ In addition, $V_{2r}$ needs to be orthogonalized and, in contrast to the other methods, the projected matrix
$S_{4r}^THS_{4r} \in \mathbb{R}^{4r \times 4r}$ needs to be computed explicitly which involves 3 matrix-matrix-multiplies of the form $W_{2r}^TXW_{2r}$ for an $n \times n$ matrix $X$.

\begin{algorithm}[ht!]
    \caption{Block $J$-orthogonal Basis Method}
    \label{alg:BJ}
    \begin{algorithmic}[1]
        \Require Hamiltonian $H{\in}\mathbb{R}^{2n\times2n}$, $u{\in}\mathbb{R}^{2n}$, $k{\in}\mathbb{N}$, $\textrm{tol}\in\mathbb{R}_{>0}$
        \Ensure  Hamiltonian $H_{4k}{\in}\mathbb{R}^{4p\times4p}$, $S_{4k}{\in}\mathbb{R}^{2n\times4p}$ with $S_{4k}^TJ_nS_{4k}=J_{2p}$, where $p \leq k$ due to the orthogonalization in line 12
        \State $u_{1} = u/\Vert u\Vert_2$
        \For {$j=2,\dots,k$}
        \State $u = Hu_{j-1}$
        \For{$i=1:j$}
            \State $z = z - ( Q(:,i)'*z)*Q(:,i)$
        \EndFor
        \If{$\Vert u\Vert_2<$ tol}
        \State stop
        \EndIf
        \State $u_{j} = u/\Vert u\Vert_2$
        \EndFor
        \State $W_{2k} = \texttt{orth}\left(\begin{bmatrix}
            \begin{bmatrix}
                (U_k)_{ij}
            \end{bmatrix}_{i,j=1}^{n} &                \begin{bmatrix}
                (U_k)_{ij}
            \end{bmatrix}_{i=n+1,j=1}^{2n}
        \end{bmatrix}\right)$
        \State $S_{4k} = \left[\begin{smallmatrix}
            W_{2k} & 0\\
            0 & W_{2k}
        \end{smallmatrix}\right]$ 
        \State $H_{4k} = S_{4k}^THS_{4k}$
    \end{algorithmic}
\end{algorithm}

\subsection{Arnoldi projection method using inner product \texorpdfstring{$x^THy$}{x'Hy}}
In \cite{LiC19}, it is noted that for symmetric positive definite $H$, it is possible to rewrite the Arnoldi method with a modified inner product with respect to $H$, i.e. $\langle \cdot,\cdot\rangle_H :=  \langle \cdot,H\cdot\rangle$.  The resulting method preserves certain first integrals to machine precision. 
As the assumption on $H$ is quite restrictive when considering the approximation of the integration of general Hamiltonian systems, we do not consider this method any further.

\subsection{Re-$J$-orthogonalization}
In theory, all of the above methods will generate a matrix  $S_r =[M~~N] \in \mathbb{R}^{2n\times 2r}$, $M=[m_1~~m_2~~\ldots~~m_r]$, $N=[n_1~~n_2~~\ldots~~n_r] \in \mathbb{R}^{2n \times r}$ with $J$-orthogonal columns. Yet, in practice, the $J$-orthogonality will be lost. Re-$J$-orthogonalization is necessary. For $j = 2, \ldots, r$, the column vectors $m_j,n_j$ of $S_r$ are re-$J$-orthogonalized by
\begin{align*}
    m_j &= m_j -S_{j-1}J_{j-1}^TS_{j-1}^TJ_nm_j,\\
    n_j &= n_j -S_{j-1}J_{j-1}^TS_{j-1}^TJ_nn_j.
\end{align*}
This re-$J$-orthogonalization is costly, it requires $16nr^2$ flops if it is applied to all column vectors of $S_r.$

\subsection{Summary of Krylov subspace methods considered}\label{subsec_summary_algorithms}
A summary of the relevant features of all methods considered is given in Table \ref{tab1}. For all methods except the block-$J$-orthogonal method, the stated computational costs refer to the generation of a (basis) matrix $S\in \mathbb{R}^{2n \times 2r}$. In the block-$J$-orthogonal method, a matrix of size  $2n \times 8r$ is computed. The cost for re(-$J$)-orthogonalization is not included. 

\begin{table}[ht!]
\caption{Summary of all Krylov subspace methods considered including computational cost$^1$ for generating a Krylov subspace of size $2r$ }
\label{tab1}      
\begin{center}
\begin{tabular}{|l|l|c|l|}
\hline
Method &  Columns  & \multicolumn{2}{c|}{Computational Cost}\\
&of $S$ are& Matvecs & Other cost \\
\hline
Arnoldi (A) & orthogonal & $2r$ & $r^2$ inner prod.\\[1ex]
Hamiltonian Lanczos (HL) & $J$-orthogonal& $2r$ & $3r$ inner prod.\\[1ex]
symplectic Arnoldi (SA) & orthogonal  & $r$& $r^2$ inner prod. \\
                        &   + $J$-orthogonal       &     & + ~$S_{2r}^THS_{2r}$\\[1ex]
isotropic Arnoldi (IA) & orthogonal & $r$ & $r^2$ inner prod. \\
                           &   + $J$-orthogonal               && \\[1ex]
Hamiltonian Extended  & $J$-orthogonal & $4r$& $14r$ inner prod. \\
Krylov (HEKS)&&&+ ~$3r$ lin. solves \\[1ex]
block $J$-orthogonal (BJ) & orthogonal & $2r$ & $r^2$ inner prod. \\
                           &   + $J$-orthogonal             &  & 3 matmats $W_{2r}^TXW_{2r}$\\
                           &                 & & + orth. $V_{2r}$\\
\hline
\end{tabular}

\footnotesize{$^1$ matvecs = matrix-vector-multiplications, matmats = matrix-matrix-multiplications}
\end{center}
\end{table}

The Arnoldi process and the block-$J$-orthogonal method 
expand the standard Krylov subspace by one vector per iteration, whereas the other symplectic Krylov subspace methods generate two or four new vectors at each step. As already outlined in the preceding sections, to ensure comparability, all methods are tuned so that a subspace of dimension $2r$ is generated. Consequently, Algorithms \ref{alg:A} and \ref{alg:BJ} are executed with $k=2r,$ while Algorithms \ref{alg:HL}, \ref{alg:SA} and \ref{alg:IA} are executed  with $k=r$.
For Algorithm \ref{alg:heks}, $k$ is chosen as $\frac{r}{2}$ in case $r$ is even and as $\frac{r-1}{2}$ otherwise. Hence, all methods generate a matrix $V_{2k} \in \mathbb{R}^{2n \times 2r}$, except for the block-$J$-orthogonal method, which produces a matrix $V_{2k} \in \mathbb{R}^{2n \times 8r}$.

Clearly, the HEKS algorithm incurs the highest computational cost, primarily due to the three linear systems that must be solved in each iteration. The block-$J$-orthogonal method and the symplectic Arnoldi method are also relatively expensive, since, in addition to the matrix–vector multiplications and inner products, further computations are required. For the remaining methods, the computational cost will depend on the sparsity pattern of $H$. If the inner products constitute the dominant part of the computational effort, then the Hamiltonian Lanczos method is likely to be the fastest, followed by the Arnoldi and the isotropic Arnoldi methods.

\section{Numerical Examples}\label{sec:numerical_experiments_expint}
The purpose of this section is to present six illustrative examples taken from the 
literature, consisting of one linear and five nonlinear partial differential equations which, after discretization, can be written as Hamiltonian systems
\begin{equation}\label{eq_HamSys}
    y'(t) = J_n^{-1}\nabla \mathcal{H}(y(t)), \qquad y(0) = y_0
\end{equation}
for a sufficiently smooth function $\mathcal{H}:\mathbb{R}^{2n} \rightarrow \mathbb{R}$, where $\nabla \mathcal{H}$ is the gradient of $\mathcal{H}.$  The numerical solution of such systems has been treated, e.g., in \cite{FQ10,HLW06,HS98,LR04,SC94}.
The Jacobian of $\nabla \mathcal{H}(y(t))$ is a symmetric matrix. Thus 
\begin{equation}\label{eq_Hx}
H_{y(t)} := J_n^{-1}D\nabla \mathcal{H}(y(t))
\end{equation}
is a Hamiltonian matrix. Solving \eqref{eq_HamSys}  numerically, using an exponential integrator,
one needs to compute $f(H)b = e^Hb$ or $f(H)b=\varphi(H)b$ for $\varphi(x)=\frac{e^x-1}{x}$.  In our study of the different Krylov subspace methods for approximating $f(H)b$, we set $H$ as the Jacobian $H_{y_0}$ of the Hamiltonian systems evaluated at 
the initial condition $y_0 \in \mathbb{R}^{2n}$. 

The linear example and two of the nonlinear test problems, a nonlinear Schrödinger equation and a nonlinear Klein–Gordon equation, are taken from \cite{EK19}. We include three further nonlinear examples adapted from \cite{MeiW17}, covering a sine–Gordon equation as well as nonlinear Schrödinger and Klein–Gordon equations. The two Schrödinger examples mainly differ in the presence of an additional linear term, a sign change, and a factor of two, while the Klein–Gordon examples vary in a few parameter choices. In all cases, the second spatial derivative is discretized using a standard central difference scheme which leads to either the discretized Laplacian 
$\Delta_n^{\text{\tiny{pbc}}}$ with periodic boundary conditions or the discretized Laplacian $\Delta_n$ with zero Dirichlet boundary conditions
\[
\Delta_n^{\text{\tiny{pbc}}} = \frac{1}{\delta_x^2}
	\begin{bmatrix}
		-2 & 1 & &1\\
		1 & \ddots & & \\
		& \ddots & \ddots & 1\\
		1& & 1 & -2
	\end{bmatrix}	,\qquad \Delta_n = \frac{1}{\delta_x^2}
	\begin{bmatrix}
		-2 & 1 & &\\
		1 & \ddots & & \\
		& \ddots & \ddots & 1\\
		& & 1 & -2
	\end{bmatrix}	.
\]
\subsection{Linear wave equation}\label{subsec:LinearWave}
The first numerical example considered is the one-dimensional wave
equation \cite[Section 6.2]{EK19}
\begin{align*}
	u_{tt}(x,t) = u_{xx}(x,t) + \frac{1}{8}\big(x(x-2)\big)^2 =  u_{xx}(x,t) + g(x),
\end{align*}
on the $x$-interval $[0,2]$ and the time interval $[0,50]$ with Dirichlet boundary conditions
$u(0,t) = u(2,t)=0,$
and initial conditions
$u(x,0) = \frac{1}{1+\sin^2(\pi x)}-1$ and $u_t(x,0) = 0.$
In order to obtain a Hamiltonian system \eqref{eq_HamSys}, a discretization is performed on an equidistant grid of size $n$ in the spatial variable $x$. The distance between two grid points is $\delta_x = 2/(n+1)$, and the grid points $x_j$ have the values  $x_j = j\delta_x$ for $j=1,\ldots,n$. Here, $n$ is chosen as $400$. Furthermore, we employ the central difference approximation to approximate the second derivative with respect to the location. This gives raise 
to the Jacobian matrix 
\begin{equation}\label{eq_Hlw}
H_{lw} = \begin{bmatrix} 0 & I_n \\ \Delta_n & 0 \end{bmatrix} \in \mathbb{R}^{2n \times 2n}.
\end{equation}
In this case, the time interval $[0,50]$ is discretized into $n_t = 2000$ subintervals,  $n_t$ time steps have to be carried out by means of numerical integration. Hence, $f(H_{lw\_ek})b$ has to be approximated for $2000$ different vectors $b$.

\subsection{Nonlinear Sine-Gordon Equation}\label{subsec:nlSGMW}
The Sine-Gordon equation  \cite[Problem 3]{MeiW17} 
\begin{align*}
	u_{tt}(x,t) &= u_{xx}(x,t)-\sin\left(u(x,t)\right),
\end{align*}
where $x\in[-5,\;5]$ and $t\in[0,\;100]$, with periodic boundary conditions
$u(-5,t) = u(5,t)$
is our second example.
This equation has been employed in a variety of physical contexts, including crystal dislocation dynamics \cite{FK39} and the propagation of magnetic flux in a Josephson junction \cite{SJ69}. Its soliton solutions, which maintain their shape during propagation and interaction, make it a prototypical example in nonlinear wave theory \cite{DJ89}.

As before, we approximate the second spatial derivative with central differences, utilizing $\delta_x=10/n$ and $x_j=-5 + j\delta_x$ for $j = 1,\ldots,n$. 
The resulting Jacobian matrix is given by
\begin{equation}\label{eq_Hsg}
H_{sg} = \begin{bmatrix}0 & I_n\\
\Delta_n^{\text{\tiny{pbc}}} + I_n& 0 \end{bmatrix}.
\end{equation}

\subsection{Nonlinear Klein-Gordon Equation}\label{subsec:nlKGE}
The two examples that follow are two variants of the nonlinear Klein-Gordon equation. The nonlinear Klein-Gordon equation is a significant tool in computational physics, with applications in nonlinear optics, quantum field theory, and solid-state physics \cite{W08}. In the latter, it is employed for the study of solitons \cite{DEM82}, a role similarly fulfilled by the sine-Gordon equation. 
Our focus will be on the cubic nonlinear Klein-Gordon equation \cite[Section 6.4]{EK19} (see also \cite[Problem 4]{MeiW17}) with periodic boundary conditions, i.e.  
\begin{align*}
\begin{split}
	u_{tt}(x,t) &= u_{xx}(x,t)-m^2u(x,t)-u^3(x,t),\\
	u(0,t) &= u(L,t),
    \end{split}
\end{align*}
where $x\in[0,\;L]$ and $t\in[0,\;T]$. 
As before, we discretize the spatial coordinate by setting $\delta_x=L/n$ and $x_j = j\delta_x$, where $j=1,\ldots,n$.

We consider this example for two different choices of the parameters $L,m$ and $n$ that can be found in the literature, leading to the  Jacobian matrix 
\begin{equation}\label{eq_Hkg}
H_{kg_\ell} = \begin{bmatrix}
0 & I_n \\ \Delta_n^{\text{\tiny{pbc}}}-B_{\ell} & 0
\end{bmatrix}, \qquad \ell \in \{1,2\}
\end{equation}
with
\[
B_1 = \frac{1}{4}I + 3\operatorname{diag}(b_1, \ldots, b_n), \quad b_j = (1+\cos(2j\pi\delta_x))^2,
\]
$n = 400$ for the first choice of parameters as in \cite[Section 6.4]{EK19}, 
and
\[
B_2 = I_n + 3\operatorname{diag}(b_1, \ldots, b_n),\quad 
b_j = (20(1+\cos(2j\pi\delta_x/1.28)))^2,
\]
$n = 512$ for the second choice of parameters as in \cite[Problem 4]{MeiW17}.

\subsection{Nonlinear Schrödinger Equation}\label{subsec:nlSGL}
The final two examples are based on the nonlinear Schrö\-din\-ger equation which describes a non-relativistic quantum-mechanical system.
The underlying equation with periodic boundary conditions  is given by 
\begin{align}\label{eq:nls_pde}
\begin{split}
    iu_t(x,t) &= \alpha u_{xx}(x,t)+\beta\vert u(x,t) \vert^2u(x,t)-\gamma\sin^2{(x)}u(x,t),\\
    u\Big(-\frac{L}{2},t\Big) &= u\Big(\frac{L}{2},t\Big),
    \end{split}
\end{align}
where $x\in [-L/2,L/2]$ and $t \in [0,T]$.
We consider this example for two different choices of the parameters that can be found in the literature. The first example is a variant of the nonlinear Schrödinger equation \cite[Example 6.3]{EK19} that describes a quasi-one-dimensional dilute gas Bose-Einstein condensate confined within a standing light wave \cite{BCDK01}.  The second variant of the nonlinear Schrödinger equation \cite[Problem 5]{MeiW17} differs in that it lacks the second linear term and features a distinct sign in the summand of the second spatial derivative.

A discretization in the spatial variable $x$ with $\delta_x = L/n$ and $x_j = -L/2+(j-1)\delta_x$, $j = 1, \ldots, n$, is used. The first example yields the Jacobian matrix 
\[
H_{ns_1} = \begin{bmatrix}
0 & -\frac{1}{2}\Delta_n^{\text{\tiny{pbc}}} -B\\ \frac{1}{2}\Delta_n^{\text{\tiny{pbc}}}+B & 0
\end{bmatrix} + \begin{bmatrix} D_2 & D_3\\-D_1&-D_2\end{bmatrix}
\]
with
\begin{align*}
    B &= \operatorname{diag}(\sin^2(x_j)), 
    &D_1 = \operatorname{diag}( 3(\hat q_0)_j^2+(\hat p_0)_j^2),\\
    D_2 &= \operatorname{diag}( 2(\hat q_0)_j (\hat p_0)_j),
    &D_3 = \operatorname{diag}( 3 (\hat p_0)_j^2+ (\hat q_0)_j^2),
\end{align*}
for \begin{align*}
	 \hat q_0 = \left[
	\operatorname{Re}\big(\sqrt{\sin^2(x_j)+1}\,e^{i\theta(x_j)}\big) \right]_{j=1}^n, \quad
	 \hat p_0 = \left[
	\operatorname{Im}\big(\sqrt{\sin^2(x_j)+1}\,e^{i\theta(x_j)}\big)\right]_{j=1}^n
\end{align*}
with $x_j = -4\pi +(j-1)\delta_x$, $\delta_x = \frac{8\pi}{500}$ and the the phase function $\theta(x)$, implicitly defined via $\tan\left(\theta\left(x\right)\right) = \pm\sqrt{2}\,\tan\left(x\right).$

For the second example we have
\[
H_{ns_2} = \begin{bmatrix}
0 & -\Delta_n^{\text{\tiny{pbc}}}\\ \Delta_n^{\text{\tiny{pbc}}} & 0
\end{bmatrix} + \begin{bmatrix} D_2 & D_3\\-D_1&-D_2\end{bmatrix}
\]
with
\begin{align*}
    D_1 = \operatorname{diag}( 6(\hat q_0)_j^2+2(\hat p_0)_j^2),\quad
    D_2 = \operatorname{diag}( 8(\hat q_0)_j (\hat p_0)_j), \quad
    D_3 = \operatorname{diag}( 6 (\hat p_0)_j^2+ 2(\hat q_0)_j^2)
\end{align*}
for 
\begin{align*} 
\begin{split}
	\hat q_0&=\left[ \operatorname{Re}\big(2e^{-i\left(2 x_j+\left(1+\frac{\pi}{2}\right)\right)}\operatorname{sech}(2x_j)\big)\right]_{j=1}^n,\\
    \hat p_0&=\left[\operatorname{Im}\big(2e^{-i\left(2 x_j+\left(1+\frac{\pi}{2}\right)\right)}\operatorname{sech}(2x_j)\big)\right]_{j=1}^n,
    \end{split}
\end{align*}
where $x_j = -10 + (j-1)\delta_x$ and $\delta_x = \frac{20}{512}$.

\section{Numerical Experiments on Approximating \texorpdfstring{$\exp(H)v$}{exp(H)v} and \texorpdfstring{$\phi(H)v$}{phi(H)v}}\label{sec5}
We compare the  Krylov subspace projection methods introduced in Section \ref{sec3} (summarized in Table \ref{tab1}) with respect to their accuracy in approximating $\varphi(hH)b$ and $e^{hH}b$, where $h\in\mathbb{R}$, using $$U_{2r}\varphi(h\widetilde{H}_{2r})W_{2r}^Tb$$ 
and $$U_{2r}e^{h\widetilde{H}_{2r}}W_{2r}^Tb,$$
respectively, where 
$$\widetilde{H}_{2r} = W_{2r}^THV_{2r}\in\mathbb{R}^{2r\times 2r},$$ 
$W_{2r}$, $U_{2r}\in\mathbb{R}^{2n\times 2r}$, $b\in\mathbb{R}^{2n}$ and $H\in\mathbb{R}^{2n\times 2n}$ with $J_nH = (J_nH)^T$. 
In case, $W_{2k} = U_{2k}$ and $U_{2k}$ has orthogonal columns or $W_{2k} = J_k^T U_{2k}^TJ_n$ and $U_{2k}$ has $J$-orthogonal columns, we have $W_{2k}^Tb = \|b\|_2e_1$ for $U_{2k}e_1 = b/\|b\|_2,$ such that the above formulae simplify.

The parameter $k$ in Algorithms \ref{alg:A} -- \ref{alg:BJ} is chosen as explained in Section \ref{subsec_summary_algorithms}. The parameter \texttt{tol} used to indicate breakdown is set to $10^{-14}$ for all algorithms and examples.

 The six different Krylov subspace methods are compared on the basis of the six examples described in the previous section, see Table \ref{tab_error_phi_phi_trick} for a short summary. For each of these examples, a Hamiltonian matrix $H=hH_j\in\mathbb{R}^{2 n_j \times 2 n_j}$, $j \in \{ lw,sg,kg_1,kg_2,ns_1,ns_2\}$ is set, where $h=0.01$ is used to mimic the situation in a timestep of an exponential integrator. The matrices $H_j$ are chosen as the Jacobian matrix $H_{y_0}= J^{-1}D\nabla \mathcal{H}(y_0)$ \eqref{eq_Hx} of one of the (non-)linear systems described in Section \ref{sec:numerical_experiments_expint}.  Recall that these matrices are indeed Hamiltonian matrices.

\begin{table}[ht!]
    \caption{Summary of all examples considered including relative error $\frac{\Vert \varphi_\text{expl}(hH)b -  \varphi_\text{impl}(hH)b\Vert_2}{\Vert\varphi_\text{impl}(hH)b\Vert_2}$ }
    \label{tab_error_phi_phi_trick} 
    \begin{center}
    \begin{tabular}{|l|l|c|c|}
        \hline
        matrix & example &  $n_j$ & rel. error  \\[0.5ex]
\hline
      $H_{lw}$&   linear wave \cite{EK19} & 400 & $5.8\cdot 10^{-13}$\\[0.5ex]
      $H_{sg}$ &    sine-Gordon equation \cite{MeiW17} & 512 & $4.2\cdot 10^{-15}$ \\[0.5ex]
       $H_{kg_1}$ &   nonlinear Klein-Gordon equation version 1 \cite{EK19} & 400& $1.2\cdot10^{-12}$\\ [0.5ex]
        $H_{kg_2}$ &   nonlinear Klein-Gordon equation version 2 \cite{MeiW17} & 512 & $7.3\cdot10^{-14}$\\ [0.5ex]
       $H_{ns_1}$ &   nonlinear Schrödinger equation version 1 \cite{EK19} & 500& $1.7\cdot 10^{-11}$\\ [0.5ex]
        $H_{ns_2}$ &  nonlinear Schrödinger equation version 2 \cite{MeiW17} & 512 & $5.6\cdot 10^{-13}$\\ [0.5ex]\hline 
    \end{tabular}
    \end{center}
\end{table}

 The vector $b\in\mathbb{R}^{2n_j}$ required for the matrix-vector products $\exp(H)b$ and $\varphi(H)b$ is chosen as a dense vector with pseudo-random entries from the standard normal distribution.  For all examples of the same size, the same vector $b$ is used to ensure consistency in the different Krylov subspaces. Clearly, this vector $b$ is also the vector used to generate each of the different Krylov subspaces. 

To compute the expression $\varphi(A)b$, we consider two variants.  The first variant involves the explicit calculation of 
\begin{align}\label{eq:varphi}
	\varphi(A)b = \big(e^A-I_{2n}\big)A^{-1}b,
\end{align}
requiring the exponential of a matrix and the solution of a system of equations. The second variant relies on the implicit calculation of the expression as observed, e.g.,  in \cite[Proposition 2.1]{S92} or \cite[Theorem 1.21]{Hig08}: for 
\begin{align*}
	\widehat{A} = \begin{bmatrix}
		A & b\\
		0& 0
	\end{bmatrix},
\end{align*}
 we obtain 
\begin{align}\label{eq:trickexp}
	e^{\widehat{A}}=\begin{bmatrix}
		e^A & \varphi(A)b\\
		0 & 1
	\end{bmatrix}.
\end{align} 
 This indicates that one can compute $\varphi(A)b$ via the exponential of a marginally larger matrix $\widehat{A}$. This approach requires minimal additional computational effort, if $e^{\widehat{A}}$ has to be computed anyway.

Table \ref{tab_error_phi_phi_trick} lists the relative difference $$\frac{\Vert \varphi_\text{expl}(A)b -  \varphi_\text{impl}(A)b\Vert_2}{\Vert\varphi_\text{impl}(A)b\Vert_2}$$ between the explicit \eqref{eq:varphi} and implicit \eqref{eq:trickexp} computation of $\varphi$ (namely $\varphi_{\text{expl}}$ and $\varphi_\text{impl}$) for the different matrices. The maximal relative difference for the examples considered is of the order of $\mathcal{O}(10^{-11})$. Therefore, we expect that both variants perform similarly for the different examples. 

 We measure the accuracy of the approximation for the six Krylov subspace methods for the approximation of $\exp(H)b$ in terms of the relative solution error 
 $$\frac{\Vert \exp(hH)b - V\exp(h\widetilde{H})W^Tb\Vert_2}{\Vert \exp(hH)b\Vert_2}.$$ 
 For the approximation of $\varphi(H)b$, we report the relative solution error 
 $$\frac{\Vert \varphi_\text{impl}(hH)b - Vf(h\widetilde{H})W^Tb\Vert_2}{\Vert \varphi_\text{impl}(hH)b\Vert_2}$$ for $f = \varphi_\text{expl}, \varphi_\text{impl}.$

All experiments are performed in MATLAB R2025b on an Intel(R) Core(TM) i7-8565U @ 1.8 GHz 1.99 GHz with 16GB RAM. The code for the numerical examples is available on Zenodo \cite{zenodo}. Our MATLAB implementation uses the standard MATLAB functions \texttt{expm} to compute the exponential of a matrix and the backslash operator \texttt{\textbackslash} to solve systems of linear equations.

\begin{figure}[ht!]
    \centering
    \includegraphics
    [width=1\textwidth]
    {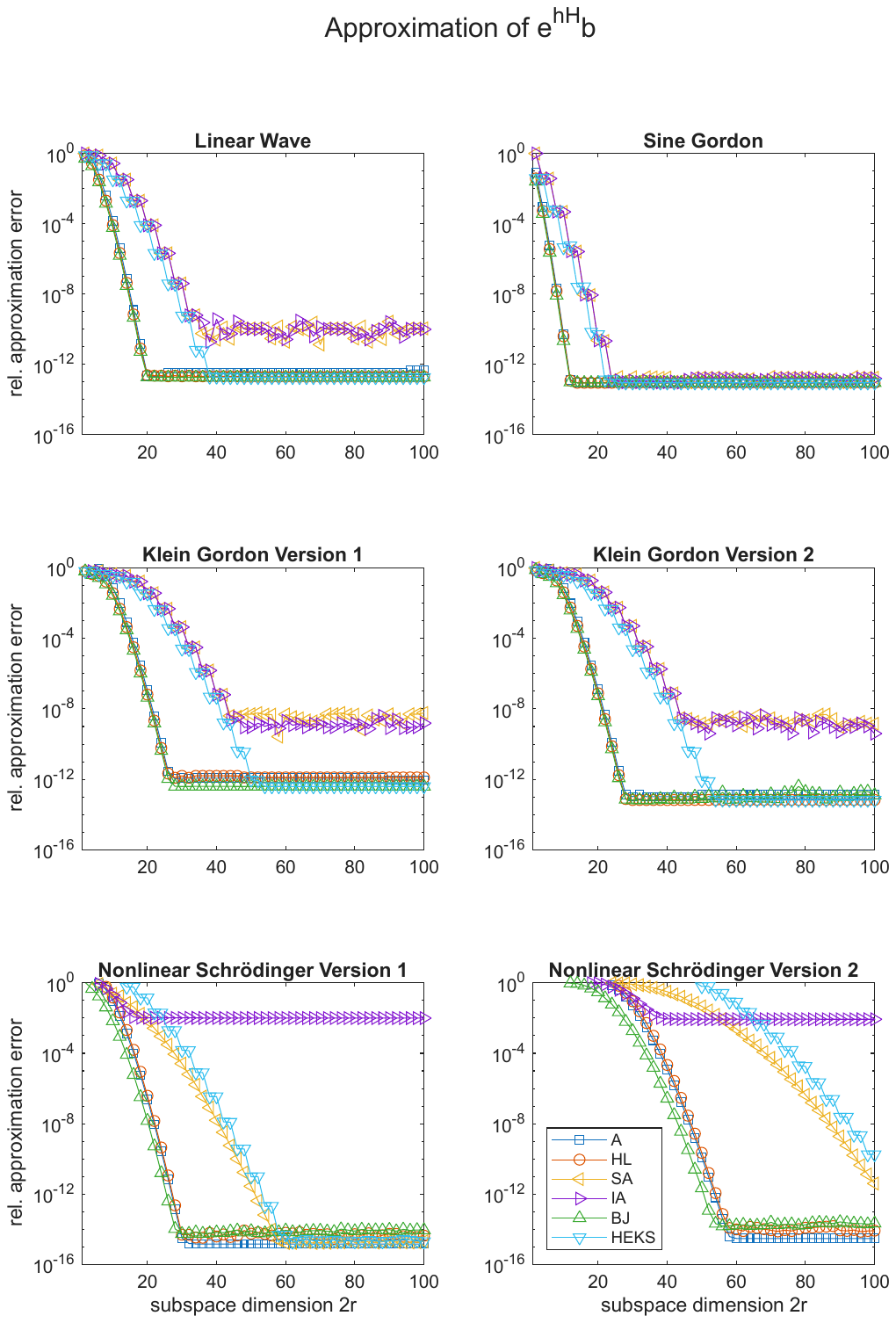}
    \caption[]{Relative solution error for different matrices in the approximation of $e^{H}b$.}
    \label{fig_approxE}
\end{figure}

\begin{figure}[ht!]
    \centering
    \includegraphics
    [width=1\textwidth]
    {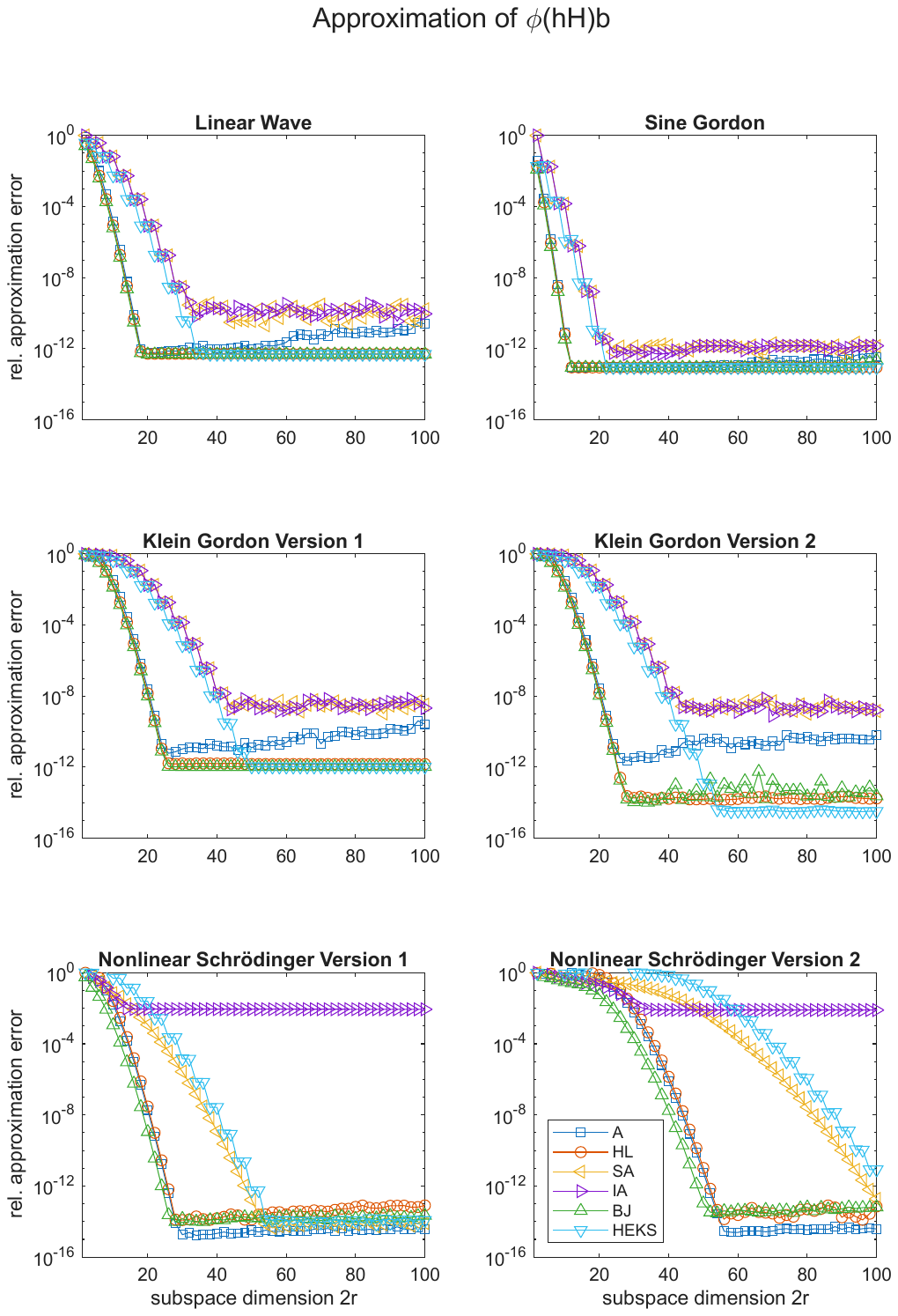}
    \caption[]{Relative solution error for different matrices in the approximation of $\varphi_\text{expl}(H)b$ using \eqref{eq:varphi}.}
    \label{fig_approxPhi}
\end{figure}

\begin{figure}[ht!]
    \centering
    \includegraphics
    [width=1\textwidth]
    {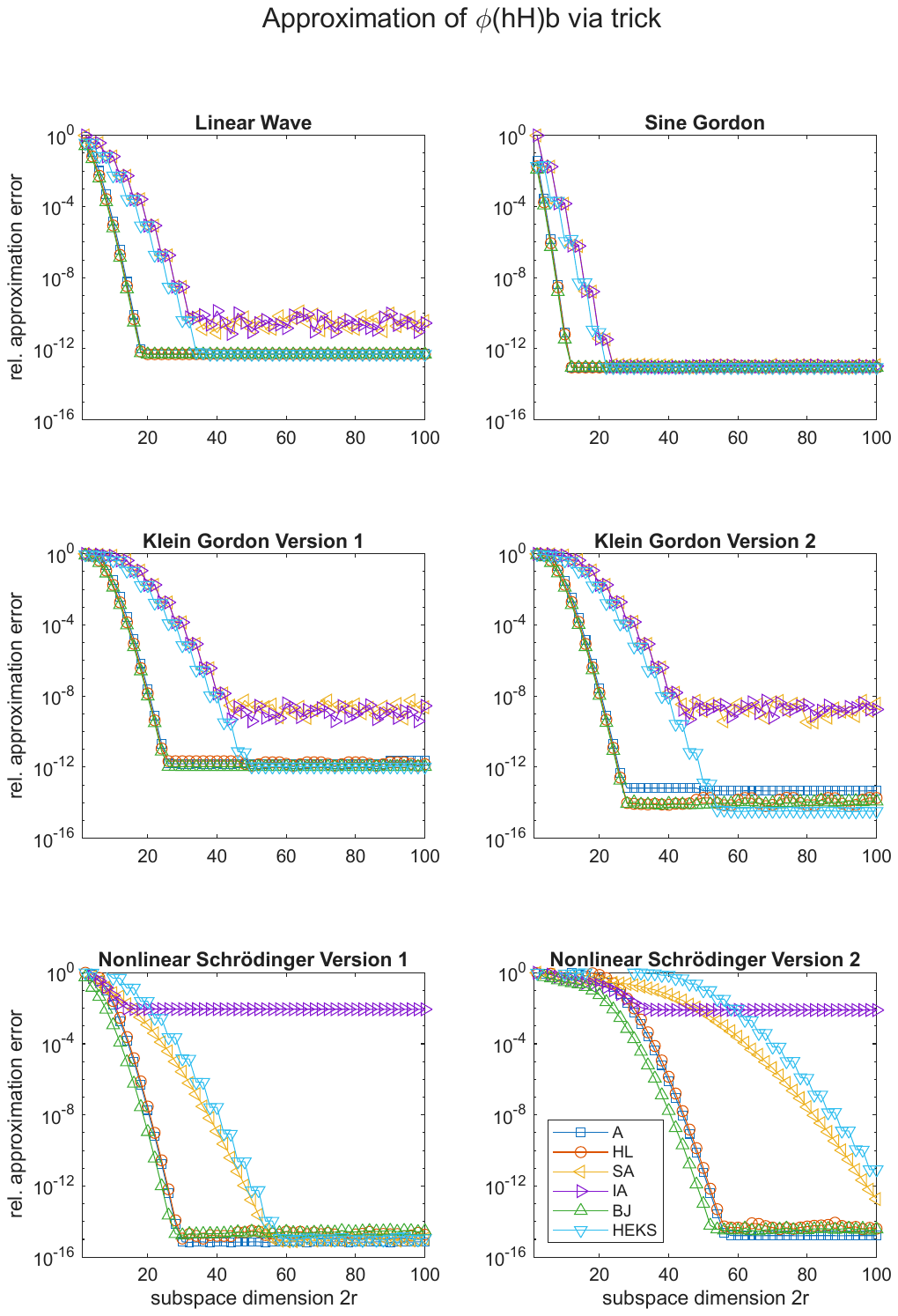}
    \caption[]{Relative solution error for different matrices in the approximation of $\varphi_\text{impl}(H)b$ using \eqref{eq:trickexp}.}
    \label{fig_approxPhiTrick}
\end{figure}

 In Fig.~\ref{fig_approxE}--\ref{fig_approxPhiTrick}, the relative solution errors of the Krylov subspace methods for approximation $e^Hb,$ $\varphi_\text{expl}(H)b$ and $\varphi_\text{impl}(H)b$ are plotted over the related subspace dimension $r = 2, 4, 6, \ldots, 50$. In practice, one should use the smallest subspace size which yields the desired accuracy. 

\begin{figure}[ht!]
    \centering
    \includegraphics
    [width=1\textwidth]
    {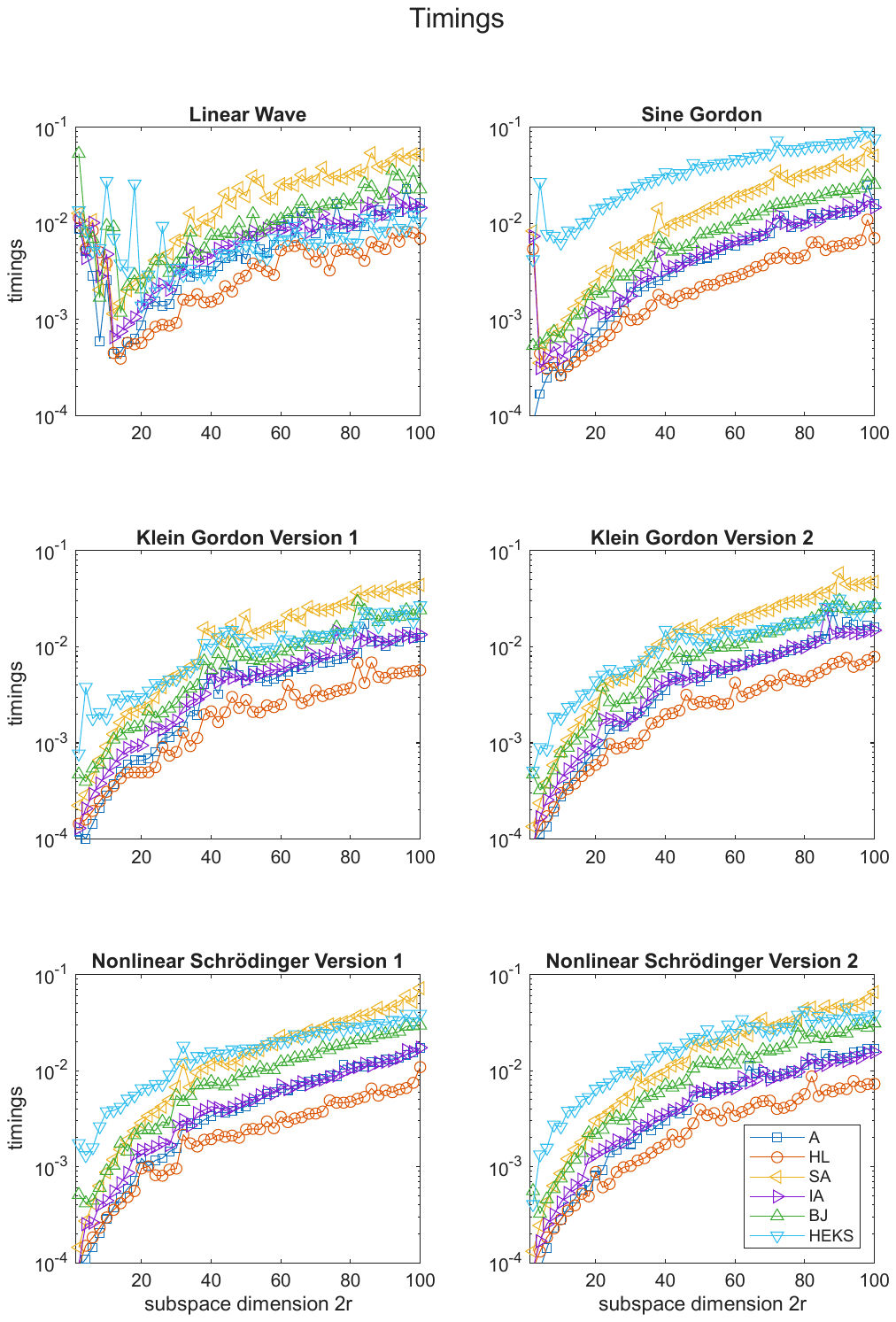}
    \caption[]{Timings for different matrices for generating the subspaces of dimension $2r$.}
    \label{fig_timings}
\end{figure}

The timings for computing the different Krylov subspaces (including re(-$J$)-orthogonalization) for the different examples can be found in Fig.~\ref{fig_timings}. The measured run-times are consistent with the findings reported in Section \ref{subsec_summary_algorithms}. Among all methods, the Hamiltonian Lanczos method achieves the shortest computation times, followed by the Arnoldi and isotropic Arnoldi methods. The Hamiltonian extended Krylov subspace method and the symplectic Arnoldi method exhibit the lowest runtime performance, whereas the block-$J$-orthogonal method 
ranks in the middle.

 As can be seen from Fig.~\ref{fig_approxE} and Fig.~\ref{fig_approxPhiTrick}, the different Krylov subspace methods behave very similarly for the approximation of $\exp(H)v$ and $\varphi_\text{impl}(H)b,$ while they show a slightly more erratic behavior for the approximation of $\varphi_\text{expl}(H)b,$ see Fig.~\ref{fig_approxPhi}.

The Arnoldi method converges to a precision of at least $\mathcal{O}(10^{-12})$ for most examples.
Among the structure-preserving methods (HL, SA, IA, HEKS, BJ), the Hamiltonian Lanczos method  and the block-$J$-orthogonal method  perform similarly to the Arnoldi method and generally converge faster than the other structure-preserving approaches. This behavior is expected, as all three methods construct a basis for the same standard Krylov subspace. The isotropic Arnoldi method  usually performs worse than the other methods, while the symplectic Arnoldi method  and the Hamiltonian extended Krylov subspace method (HEKS) typically converge much more slowly than the Hamiltonian Lanczos, the block-$J$-orthogonal and the Arnoldi method.
All methods except the Arnoldi and Hamiltonian Lanczos methods contain directions in their bases that extend beyond the standard Krylov subspace. Based on findings from Fig.~\ref{fig_approxE}--\ref{fig_approxPhiTrick}, it seems reasonable to suggest that a method for approximating the functions in question solely considering the standard Krylov subspace is a suitable approach. 

Considering both convergence behavior and computational cost, the Hamiltonian Lanczos method is the most efficient (structure-preserving) approach, followed by the numerically stable Arnoldi method. Although the Hamiltonian Lanczos method is inherently unstable, this instability can be monitored and controlled during the computation and, therefore, does not cause practical difficulties.

\clearpage

\subsection{Adaptive selection of the Krylov subspace dimension}
As noted previously, in practical applications one does not necessarily work with a fixed subspace dimension; instead, it may be chosen adaptively depending on the problem. The following discussion outlines how this can be realized in our implementation. 

In the literature, particularly in the context of numerical exponential integrators, numerous approaches have been proposed that combine approximation $\varphi(A)b$ and $e^Ab$ by building a Krylov subspace using the Arnoldi iteration, together with time-stepping schemes, in order to prevent the Krylov subspace from growing excessively. These fully adaptive methods adjust both the time-step size and the dimension of the Krylov subspace to achieve the desired accuracy; see, e.g., the works
on adaptive selection of the Krylov subspace size \cite{NW11,NW12,LTT12,GRT18,BerS24}.
For the purposes of this paper, these approaches are not suitable. Instead, a rather simple adaptive strategy based on an error estimator $\| e^H b - U_k e^{\widehat{H}_k}e_1\|$ or $\| \varphi(H) b - U_k \varphi({\widehat{H}_k})e_1\|$ as in \cite{S92} can be employed.
Ideas on how to use restarted Krylov methods for approximating a matrix function times a vector can be found, e.g., in \cite{EE06,AEEG08,EEG11,FGS14}. 

In the context of approximating $e^Hb$ through the Arnoldi method, the approximation error can be specified exactly.
\begin{theorem}[{\cite[Theorem 5.1]{S92}}] \label{theo:error}
	The error produced by the Arnoldi  approximation \eqref{eq:FuncApprox} satisfies the following expansion:
	\begin{align}\label{eq:FuncApproxError}
		e^Hu_1 - U_k e^{\widehat{H}_k}e_1=\widehat h_{k+1,k}\sum_{i=1}^{\infty}e_k^T\varphi_k(\widehat{H}_i)e_1H^{i-1}u_{k+1}
	\end{align}
    for $\varphi_0(z) = e^z$ and $\varphi_{i+1}(z) = \frac{\varphi_i(z)-\varphi_i(0)}{z}, i \geq 1$, and $u_1 = b/\|b\|_2$ using the notation as in \eqref{eq:arnoldi_recursion}.
\end{theorem}
As noted in \cite{S92}, the first one or two summands of the sum in \eqref{eq:FuncApproxError} are sufficient to provide an appropriate estimate of the error \eqref{eq:FuncApproxError}. Thus, for the approximation via the Arnoldi method, we employ the value of the first summand 
\begin{align}\label{err_A}
	\epsilon_k^\texttt{A} = \Vert b\Vert_2 \vert h\widehat{h}_{k+1,k}e_k^T\varphi(h\widehat{H}_k)e_1\vert,
\end{align}
as the error estimator for the approximation of $e^{hH}b$. In \cite[Section 6.3]{HLS98}, the same error estimate, but for the $\varphi$-function, is derived using an argumentation via the Cauchy integration rule.

An inspection of the proof of Theorem \ref{theo:error} in \cite{S92} shows that is also applicable to the Hamiltonian Lanczos method. Hence, the two algorithms that performed best in terms of accuracy and computational efficiency can both be implemented in an adaptive manner. Using the notation as in \eqref{eq:hamlan_recursion}, we use
\begin{align}\label{err_HL}
	\epsilon_k^\texttt{HL} = \Vert b\Vert_2 \vert h\beta_{k}e_{2k}^T\varphi(h{H}_{2k})e_1\vert
\end{align}
as the error estimator for the approximation of $e^{hH}b$ by the Hamiltonian Lanczos method. 

Replacing the for-loop over $j$ in Algorithms \ref{alg:A} and \ref{alg:HL} with a while-loop that terminates once \eqref{err_A} or \eqref{err_HL}, respectively, falls below a prescribed tolerance yields an adaptive version of the algorithms. Table \ref{tab_adaptive} presents the exact errors and the estimated errors given by \eqref{err_A} and \eqref{err_HL} for the approximation of $e^Hb$ in the test problems $H_{kg_1}$ and $H_{ns_2}$. As already observed in \cite{S92}, the error estimate for the approximation via the Arnoldi method are surprisingly sharp. The same holds true for the error estimate for the approximation via the Hamiltonian Lanczos method.

\begin{table}[ht!]
    \caption{Actual error and estimate for $H_{kg_1}$ and $H_{ns_2}$}\label{tab_adaptive}
    \begin{center}
    \begin{tabular}{|c|c|c|c|c|}
    \hline
    \multicolumn{5}{|c|}{Klein-Gordon $H_{kg_1}$}\\ \hline
    & \multicolumn{2}{c|}{Arnoldi method} & \multicolumn{2}{c|}{Hamiltonian Lanczos method}\\
    k & Actual error & $\epsilon_k^\texttt{A}$ &Actual Error & $\epsilon_k^\texttt{HL}$\\ \hline  
 $1$&$  1.2410\cdot 10^{+01} $&$  1.1376\cdot 10^{-01} $&$ 5.8950\cdot 10^{+00} $&$  9.1674\cdot 10^{+01}$\\
 $2$&$  8.4214\cdot 10^{-02} $&$  1.4287\cdot 10^{-03} $&$ 4.9892\cdot 10^{-02} $&$  1.3203\cdot 10^{-01}$\\
 $3$&$  3.5103\cdot 10^{-04} $&$  8.0418\cdot 10^{-06} $&$ 1.8852\cdot 10^{-04} $&$  1.7009\cdot 10^{-04}$\\
 $4$&$  7.8315\cdot 10^{-07}  $&$ 2.2719\cdot 10^{-08} $&$ 4.2073\cdot 10^{-07} $&$  7.0821\cdot 10^{-07}$\\
 $5$&$  1.0913\cdot 10^{-09} $&$  3.9726\cdot 10^{-11} $&$ 6.2677\cdot 10^{-10} $&$  5.8256\cdot 10^{-10}$\\
 $6$&$  6.6843\cdot 10^{-11} $&$  4.8277\cdot 10^{-14} $&$ 6.7101\cdot 10^{-11} $&$  8.5729\cdot 10^{-13}$\\ \hline
    \end{tabular}
    \bigskip
    
        \begin{tabular}{|c|c|c|c|c|}
        \hline
    \multicolumn{5}{|c|}{Nonlinear Schrödinger $H_{ns_2}$}\\ \hline
   & \multicolumn{2}{c|}{Arnoldi method} & \multicolumn{2}{c|}{Hamiltonian Lanczos method}\\
    k & Actual error & $\epsilon_k^\texttt{A}$ &Actual Error & $\epsilon_k^\texttt{HL}$\\ \hline  
$1$&$    9.2721\cdot 10^{-01}$&$   1.0300\cdot 10^{+00}$&$    1.0537\cdot 10^{+00}$&$   1.4911\cdot 10^{+00}$\\
$2$&$   1.3030\cdot 10^{-01}$&$   1.3712\cdot 10^{-01}$&$  1.8027\cdot 10^{-01}$&$   3.1835\cdot 10^{-01}$\\
$3$&$   7.8115\cdot 10^{-03}$&$   8.0209\cdot 10^{-03}$&$ 1.1104\cdot 10^{-02}$&$   2.2716\cdot 10^{-02}$\\
$4$&$   2.4533\cdot 10^{-04}$&$   2.4964\cdot 10^{-04}$&$ 3.6175\cdot 10^{-04}$&$   6.9830\cdot 10^{-04}$\\
$5$&$   4.9355\cdot 10^{-06}$&$   4.9948\cdot 10^{-06}$&$7.6003\cdot 10^{-06}$&$   1.4510\cdot 10^{-05}$\\
$6$&$   6.4087\cdot 10^{-08}$&$   6.4648\cdot 10^{-08}$&$  1.0149\cdot 10^{-07}$&$   1.9719\cdot 10^{-07}$\\
$7$&$   6.0482\cdot 10^{-10}$&$   6.0899\cdot 10^{-10}$&$    9.6888\cdot 10^{-10}$&$   1.9068\cdot 10^{-09}$\\
$8$&$   4.4909\cdot 10^{-12}$&$   4.5155\cdot 10^{-12}$&$  7.2386\cdot 10^{-12}$&$   1.4340\cdot 10^{-11}$\\
$9$&$   2.5424\cdot 10^{-14}$&$   2.5258\cdot 10^{-14}$&$   4.0502\cdot 10^{-14}$&$   8.3713\cdot 10^{-14}$\\ \hline
    \end{tabular}
    \end{center}
\end{table}

\section{Conclusions}\label{sec6}
In this work, we have conducted a comparative study of several structure-preserving Krylov subspace methods for approximating the matrix–vector products \eqref{eq_problem}, where 
$H$ is a large and sparse Hamiltonian matrix. Our focus was on methods that preserve the Hamiltonian/symplectic structure of the problem, since this property is crucial in the context of exponential integrators for Hamiltonian systems.

The numerical experiments, performed on a representative set of Hamiltonian test problems, demonstrate clear differences in performance between the considered approaches. Among all structure-preserving methods, the Hamiltonian Lanczos algorithm showed the most favorable balance between computational efficiency and accuracy. In particular, it achieved high-quality approximations with significantly reduced computational time comparable to the non-structure-preserving Arnoldi method. Furthermore, we have illustrated how the dimension of the Krylov subspace can be chosen adaptively for the Arnoldi and the Hamiltonian Lanczos algorithm, allowing for a reliable accuracy control without excessive computational effort. The availability of such an adaptive variant of the Hamiltonian Lanczos method further enhances its practical applicability in large-scale problems. In summary, the non-structure-preserving Arnoldi method and the structure-preserving Hamiltonian Lanczos method emerged as the two most promising approaches in our comparative study.


\section*{Declarations}

\begin{itemize}
\item Conflict of interest/Competing interests: -none-
\item Consent for publication: All three authors have read the final version of the manuscript and have given their consent for publication.
\item Code availability: see \cite{zenodo} 
\item Author contribution; Peter Benner: Conceptualization, Writing – review \& editing; Heike Faßbender: Conceptualization, Methodology, Supervision, Writing - Original Draft;  Michel-Niklas Senn: Methodology, Software, Investigation, Validation, Writing - Original Draft
\item Originality/Prior Publication: This work is based in parts on \cite{Sen25}.
\end{itemize}


\end{document}